\input amstex
\input amsppt.sty
\TagsOnRight

\NoBlackBoxes

\define\Y{\Bbb Y}
\define\Z{\Bbb Z}
\define\C{\Bbb C}
\define\R{\Bbb R}

\define\de{\delta}

\define\la{\lambda}
\define\si{\sigma}
\define\th{\theta}

\define\X{\frak X}
\define\XX{\wt X}

\define\wt{\widetilde}

\define\tht{\thetag}

\define\Prob{\operatorname{Prob}}

\define\const{\operatorname{const}}

\define\Conf{\operatorname{Conf}}

\define\kr{{\text{Krawtchouk}}}
\define\ch{{\text{Charlier}}}
\define\he{{\text{Hermite}}}

\define\dsine{{\text{discrete sine}}}
\define\Mix{\text{Mix}}

\define\Dim{\operatorname{Dim}}
\define\SW{{\text{Schur-Weyl}}}
\define\PSW{{\text{Poisson-Schur-Weyl}}}
\define\Plancherel{{\text{Plancherel}}}

\topmatter
\title Asymptotics of Plancherel--type random partitions \endtitle

\author Alexei Borodin and Grigori Olshanski \endauthor

\dedicatory Dedicated to E.~B.~Vinberg \enddedicatory

\abstract We present a solution to a problem suggested by Philippe Biane: We
prove that a certain Plancherel--type probability distribution on partitions
converges, as partitions get large, to a new determinantal random point process
on the set $\Z_+$ of nonnegative integers. This can be viewed as an edge limit
transition. The limit process is determined by a correlation kernel on $\Z_+$
which is expressed through the Hermite polynomials, we call it the discrete
Hermite kernel. The proof is based on a simple argument which derives
convergence of correlation kernels from convergence of unbounded self--adjoint
difference operators.

Our approach can also be applied to a number of other probabilistic models. As
an example, we discuss a bulk limit for one more Plancherel--type  model of
random partitions.
\endabstract

\endtopmatter

\document

\head Introduction \endhead

This work appeared as our attempt to solve a problem posed by Philippe Biane.
In \cite{Bi2} he considered a model of random partitions arising from
decomposition of tensor spaces $(\C^N)^{\otimes n}$ (the Schur--Weyl duality
between representations of the symmetric group $S_n$ and the unitary group
$U(N)$). The partitions in question have at most $N$ nonzero parts which sum to
$n$, and the weight of a partition $\la$ is proportional to the product of
dimensions of the irreducible representations of $S_n$ and $U(N)$ indexed by
$\la$.

Biane discovered that as $n$ and $N$ go to infinity so that $\sqrt n\sim c N$
then the boundary of the Young diagram associated to the random partition
$\la$, suitably scaled, tends to a nonrandom limit curve given by an explicit
formula. The limit curve depends on  the parameter $c>0$.

If $n$ is fixed while $N\to\infty$ then the model turns into the well--known
Plancherel model of random partitions of $n$. This agrees with the fact that as
$c$ approaches 0, Biane's limit curve turns into the celebrated
Vershik--Kerov--Logan--Shepp limit curve for the Plancherel model found in
\cite{VK1}, \cite{VK2}, \cite{LS}.

Biane's formulas show that the value $c=1$ is special: The tangent line to the
limit curve at one of its endpoints has (in appropriate coordinates) slope $-1$
for $c<1$, 0 for $c=1$, and $+1$ for $c>1$. Biane's question concerned the
local structure of the boundary of the random Young diagram at $c=1$ near this
point of the limit shape.

We address this question in a modified form. Namely, we replace the initial
probability distribution on partitions by its poissonization with respect to
parameter $n$. This procedure is well known, we explain it in \S4. One expects
that poissonization does not affect the asymptotics but we leave the discussion
of this issue out of this paper.

After poissonization, the probability distribution lives on all partitions
$\la=(\la_1,\la_2,\dots,\la_N)$ with at most $N$ nonzero parts and without any
constraints on $|\la|=\la_1+\dots+\la_N$. It is convenient to interpret $\la$
as an $N$--point configuration on $\Z_+=\{0,1,2,\dots\}$ via
$$
\la\mapsto\{x_1,\dots,x_N\}, \qquad x_i=\la_i+N-i.
$$
The weight of $\la$ now depends on the poissonization parameter $\nu>0$ (which
replaces $n$) and has the form
$$
\const\cdot \prod_{i=1}^N\frac{(\nu/N)^{x_i}}{x_i!}\cdot \prod_{1\le i<j\le N}
(x_i-x_j)^2.
$$
This is the so--called Charlier orthogonal polynomial ensemble.

Our main result is the following statement (see below Theorem 4.1 and
Proposition 3.3).

\proclaim{Theorem} Fix an arbitrary $s\in\R$. Let $N=1,2,\dots$ and assume that
the parameter $\nu=\nu(N)$ depends on $N$ in such a way that
$$
\nu=N^2+(s+o(1))N^{3/2}, \qquad N\to\infty.
$$

As $N\to\infty$, the probability distribution of $\{x_1,\dots,x_N\}$ converges
to a probability measure on $2^{\Z_+}$, the set of all point configurations $Y$
on $\Z_+$. The correlation functions of the limit measure have the form
{\rm(}$k=1,2,\dots${\rm)}
$$
\Prob\{Y\supset \{y_1,\dots,y_k\}\}=\det[K_s(y_i,y_j)]_{1\le i,j\le k}\,,
\qquad y_1,\dots,y_k\in\Z_+,
$$
where, for $x,y\in\Z_+$,
$$
\multline K_s(x,y)=(\pi x!y!2^{x+y})^{-\frac12}\,
\int_{s/\sqrt2}^{+\infty}e^{-t^2}H_x(t)H_y(t)dt\\
=(\pi
x!y!2^{x+y})^{-\frac12}\,e^{-\frac{s^2}2}\,\,\frac{xH_{x-1}(\tfrac{s}{\sqrt2})\cdot
H_y(\tfrac{s}{\sqrt2}) -H_x(\tfrac{s}{\sqrt2})\cdot
yH_{y-1}(\tfrac{s}{\sqrt2})}{x-y}\\
=(4\pi
x!y!2^{x+y})^{-\frac12}\,e^{-\frac{s^2}2}\,\,\frac{H_{x+1}(\tfrac{s}{\sqrt2})
H_y(\tfrac{s}{\sqrt2}) -H_x(\tfrac{s}{\sqrt2})
H_{y+1}(\tfrac{s}{\sqrt2})}{x-y}\,.
\endmultline
$$
Here $H_m$ is the classical Hermite polynomial, see \cite{KS}.
\endproclaim

The introduction of the additional parameter $s$ above is also due to Biane.

The determinantal structure of the correlation functions means that the limit
measure belongs to the class of determinantal random point processes which
arise in a variety of probabilistic models, see, e.g., \cite{So1}, \cite{So2},
\cite{BHKPV}, \cite{Ly}. In particular, the determinantal processes arise in
connection with the Plancherel measure, see \cite{Jo1}, \cite{BOO}.

We call $K_s(x,y)$ (the kernel of the determinantal formula) the {\it discrete
Hermite kernel\/}. Many similar examples of correlation kernels are known,
however, to our best knowledge, the discrete Hermite kernel is new.

To prove the theorem we have to check that the Charlier correlation kernel
(which is essentially the $N$th Christoffel--Darboux kernel for the orthogonal
Charlier polynomials) converges to the discrete Hermite kernel. Usually such
facts are verified using asymptotics of orthogonal polynomials (see, e.g.,
\cite{Jo1} for a different limit regime for the Charlier kernel). Such an
approach is applicable to our problem as well. However, we take another path
and extract the needed convergence from an abstract theorem concerning strong
resolvent convergence of unbounded self--adjoint operators. These self--adjoint
operators appear as difference operators on $\Z_+$ associated to the Charlier
polynomials. This approach seems to be new \footnote{Even though the method of
deriving the asymptotics of special functions through differential equations is
well known, we have never seen the same idea applied to correlation kernels.},
and it appears to be much less technical comparing to the traditional one.

To demonstrate the effectiveness of this approach we apply it to another model
of representation--theoretic origin. This model appeared in the works of Biane
\cite{Bi1} and of Pittel and Romik \cite{PR}, it turns out to be related to the
so--called Krawtchouk orthogonal polynomial ensemble. In \S5 we sketch a proof
of the convergence of the Krawtchouk kernel to the discrete sine kernel. This
result cannot be viewed as new: it can be extracted from \cite{IS},
\cite{BKMM}. The point is that our argument is short and direct. We show that
the result allows one to predict the form of the limit shape obtained in
\cite{Bi1} and \cite{PR}.

\subhead Acknowledgement \endsubhead We are very grateful to Philippe Biane for
posing the problem and sharing his insights with us. We are also grateful to
Barry Simon for helpful advice. The research of A.~Borodin was partially
supported by the NSF grant DMS-0402047. Both authors were also supported by the
CRDF grant RUM1-2622-ST-04.

\head 1. Preliminaries and the Plancherel model
\endhead

Let $\Z_+=\{0,1,2,\dots\}$ be the set of nonnegative integers. Recall that a
{\it partition\/} is an infinite sequence $\la=(\la_1,\la_2,\dots)$ of
nonincreasing numbers from $\Z_+$ with finitely many nonzero terms. The sum of
the terms is denoted as $|\la|$. We say that $\la$ is a partition of $n$ if
$n=|\la|$.

Following \cite{Ma} we identify partitions and Young diagrams. A {\it Young
diagram\/} is a finite collection of unit squares in the quarter plane with
coordinates $(i,j)$, where the $i$--axis is directed downward and the $j$--axis
is directed to the right; a square with lower right corner
$(i,j)\in\{1,2,\dots\}\times\{1,2,\dots\}$ enters the diagram of $\la$ if and
only if $\la_i\ge j$. Thus, the diagram $\la$ has $\la_i$ squares in the $i$th
row with the numbering of rows ranging from the top to the bottom, and the
total number of squares is equal to $|\la|$. The set of all partitions (=Young
diagrams) will be denoted as $\Y$ and the subset of partitions of $n\in\Z_+$
will be denoted as $\Y_n$.

The {\it conjugate\/} partition $\la'$ is obtained, in terms of Young diagrams,
by transposing the coordinate axes. Clearly, $\la'_1$ coincides with the number
of nonzero terms $\la_i$ (or the number of nonvoid rows in the diagram); this
number is also denoted as $\ell(\la)$.

The {\it boundary\/} of a diagram $\la$ is a broken line going from the point
$(i,j)=(\la'_1,0)$ to the point $(i,j)=(0,\la_1)$. It is convenient to add to
the boundary those parts of the coordinate axes that are below $(\la'_1,0)$ and
to the right of $(0,\la_1)$. The boundary of $\la$ will be denoted as
$\partial\la$.

Partitions $\la\in\Y$ can also be regarded as {\it particle configurations\/}
on a 1--di\-men\-sional lattice. Let $\Z'=\Z+\frac12$ denote the set of
(proper) half--integers. We assign to $\la$ an infinite subset of $\Z'$:
$$
X(\la)=\{\la_i-i+\tfrac12\}, \qquad i=1,2,\dots, \tag1.1
$$
and we regard $X(\la)$ as a configuration of {\it particles\/} sitting at nodes
of the lattice $\Z'$. The unoccupied nodes of $\Z'$ will be called {\it
holes\/}. Note a duality relation between particles and holes: reflecting the
configuration of holes about 0 we get $X(\la')$.

The boundary $\partial\la$ (with parts of the coordinate axes included) can be
viewed as a doubly infinite sequence of horizontal and vertical unit segments,
and the particle configuration $X(\la)$ is a convenient way to encode that
sequence. Specifically, a node $k\in\Z'$ is occupied by a particle from
$X(\la)$ if and only if the the line $j-i=k$ meets the boundary $\partial\la$
at the midpoint of a vertical segment. Likewise, the holes correspond to the
midpoints of horizontal segments. This correspondence makes evident the
particle/\-hole duality mentioned above.

Let $S_n$ denote the symmetric group of degree $n$. The irreducible
$S_n$--modules are parametrized by the Young diagrams with $n$ squares. Recall
the set of such diagrams is denoted as $\Y_n$. For an arbitrary diagram
$\la\in\Y_n$, let $\dim\la$ denote the dimension of the corresponding
irreducible module. Equivalently, $\dim\la$ equals the number of standard Young
tableaux of the shape $\la$.

By Burnside's theorem,
$$
\sum_{\la\in\Y_n}(\dim\la)^2=n!
$$
The {\it Plancherel measure\/} on $\Y_n$, denoted as $M^\Plancherel_n$, is
defined as the probability measure with the weights
$$
M^\Plancherel_n(\la)=\frac{(\dim\la)^2}{n!}\,, \qquad \la\in\Y_n\,.
$$

Regard diagrams $\la\in\Y_n$ as random objects defined on the probability space
$(\Y_n,M^\Plancherel_n)$. As $n\to\infty$, the boundary of the random diagram,
scaled by the factor of $n^{-1/2}$, tends to a (nonrandom) limit curve. This is
a well--known result due to Logan--Shepp \cite{LS} and Vershik--Kerov
\cite{VK1}, \cite{VK2} (see also Kerov's book \cite{Ke2}). Specifically, if $i$
and $j$ denote the initial row and column coordinates then the scaled
coordinates are defined as $x=i\cdot n^{-1/2}$ and $y=j\cdot n^{-1/2}$, and the
equation of the limit curve in the $(x,y)$ plane can be written as
$$
x+y=\Omega(y-x), \qquad -2\le y-x\le2,
$$
where
$$
\Omega(u)=\tfrac2\pi(u\arcsin\tfrac u2+\sqrt{4-u^2}).
$$

This result leads to the following conclusions:

(a) Observe that the limit curve meets the coordinate axes $x=0$ and $y=0$ at
points $(0,2)$ and $(2,0)$ respectively. This suggests that the first row
$\la_1$ and the the first column $\la'_1$ of the typical Plancherel diagram
$\la\in\Y_n$ should grow as $2\sqrt n$, which is indeed true: see \cite{VK1},
\cite{VK2} for a precise statement. Moreover, the same holds for each of the
largest row and column lengths $\la_k$, $\la'_k$, where the index $k$ is
arbitrary but fixed.

(b) Let, as above, $i$ and $j$ be the row and column coordinates.  Fix
$u\in(-2,2)$ and let $a_n$ be a sequence of positive numbers such that
$a_n\to\infty$ and $a_n/\sqrt n\to0$. Then, as $n$ gets large, the proportion
of horizontal (respectively, vertical) steps of the random boundary, contained
in the strip $|j-i-u\cdot\sqrt n|\le a_n$, should be close to
$(1+\Omega'(u))/2$ (respectively, to $(1-\Omega'(u))/2$), where $\Omega'(u)$ is
the derivative of $\Omega(u)$. This statement just means that the slope of the
boundary of our random Young diagram approximates the slope of the limit curve.

A somewhat different but essentially equivalent formulation is as follows: Let
$k_n$ be a sequence of half--integers such that $k_n/\sqrt n\to u\in(-2,2)$.
Given $n$, look at the intersection of the line $j-i=k_n$ with the boundary of
the random diagram $\la\in\Y_n$. This is a midpoint of a boundary segment,
which can be either horizontal or vertical. Then, for $n$ large, the
probability to find a horizontal segment should be close to $(1+\Omega'(u))/2$.
This is indeed true, see \cite{BOO}.

(c) As $u$ ranges over $(-2,2)$, the quantity $(1+\Omega'(u))/2$ monotonically
increases, so that that the probability of finding horizontal fragments
increases, too. At the endpoints $-2$ and $2$ the quantity $(1+\Omega'(u))/2$
takes values 0 and 1 (that is, the limit curve is tangent to the coordinate
axes $y=0$ and $x=0$). This suggests that, typically, each of the differences
$\la_k-\la_{k+1}$ or $\la'_k-\la'_{k+1}$ (with $k$ fixed) should grow as
$n\to\infty$. A much more precise statement can be found in \cite{Ok},
\cite{BOO}, \cite{Jo1}. In particular, it turns out that the order of growth of
these differences is $n^{1/6}$.

In the next sections we will consider two other models of random Young diagrams
which may be viewed as deformations of the Plancherel model.

\head 2. Biane's model \endhead

There is a close relationship between the Plancherel model and the biregular
representation of the symmetric group. Indeed, the group $S_n$ acts on itself
by left and right shifts. The corresponding representation of the group
$S_n\times S_n$ in the space of functions on $S_n$  has simple spectrum indexed
by diagrams $\la\in\Y_n$, and $(\dim\la)^2$ is just the dimension of the
irreducible component indexed by $\la$. Since $n!$ is the dimension of the
whole representation space, the Plancherel weight of $\la$ is equal to the
relative dimension of the irreducible component indexed by $\la$.

Now we apply the same construction but starting with a different representation
with simple spectrum. Let $n$ and $N$ be two positive integers. Consider the
tensor space $(\C^N)^{\otimes n}$ as a bimodule with respect to the natural
commuting actions of the symmetric group $S_n$ and the unitary group $U(N)$. By
the Schur--Weyl duality, the representation of the group $S_n\times U(N)$ in
$(\C^N)^{\otimes n}$ has simple spectrum which is indexed by Young diagrams
$\la$ such that $|\la|=n$ and $\ell(\la)\le N$. Let $\Y_n(N)$ stand for the set
of such diagrams. For $\la\in\Y_n(N)$, the dimension of the corresponding
irreducible component equals $\dim\la\cdot\Dim_N\la$, where by $\Dim_N\la$ we
denote the dimension of the irreducible polynomial $U(N)$--module with highest
weight $(\la_1,\dots,\la_N)$. This serves as a prompt for introducing a
probability measure on $\Y_n(N)$:
$$
M^\SW_{n,N}(\la)=\frac{\dim\la\cdot\Dim_N\la}{N^n}\,, \qquad \la\in\Y_n(N)
\tag2.1
$$
(the factor $N^n$ in the denominator is the dimension of the whole tensor
space).

Let us take $M^\SW_{n,N}$ as the distribution law for a random ensemble of
diagrams $\la\in\Y_n(N)$ and ask about the asymptotic properties of this
ensemble as $n$ and $N$ go to infinity.

For the first time, this question was addressed by Sergei Kerov \cite{Ke1} (see
also \cite{Ke2, Chapter III, \S3}). He showed that if $n$ and $N$ have the same
order of growth (that is, $n/N$ tends to a positive constant) then, after
scaling by the factor of $n^{-1/2}$, the boundary of the random diagram tends
to a limit shape, which is exactly the same as in the case of the Plancherel
model. This result admits the following heuristic explanation:

If $N\ge n$ then $\Y_n(N)$ coincides with $\Y_n$, and it is readily checked
that
$$
\lim_{N\to\infty}M^\SW_{n,N}(\la)\to M^\Plancherel_n(\la), \qquad \la\in\Y_n\,.
$$
On the other hand, a typical Plancherel diagram $\la\in\Y_n$ has approximately
$2\sqrt{n}$ rows, which explains why the constraint of type $\ell(\la)\le
N=O(n)$ turns out to be asymptotically negligible.

Finer results were obtained later by Philippe Biane \cite{Bi2}. \footnote{We
strongly encourage the reader to look at this paper for a better understanding
of what follows.} He examined a family of limit regimes depending on parameter
$c\in(0,+\infty)$:
$$
n\to\infty, \qquad N\sim c^{-1}\sqrt n,
$$
and discovered that, for $c$ fixed, the scaled random diagrams concentrate near
a limit curve $x+y=P_c(y-x)$ depending on $c$. The curves $v=P_c(u)$ are
explicitly described in \cite{Bi2, \S3.1}, they provide an interesting
deformation of the Plancherel curve $v=\Omega(u)$, which appears as the limit
case $c=0$.

Look at the intersection of the curve $v=P_c(u)$ with the line $v=-u$, which
happens at $u=c-2$. A close examination of Biane's formulas (see the end of
\S3.1 in \cite{Bi2}) reveals the following fact:
$$
\left.\frac{dP_c(u)}{du}\right|_{u=c-2}=\cases -1, & c<1\\0, & c=1\\ +1, & c>1
\endcases.
$$

Our interest is what happens at the critical value $c=1$ corresponding to the
limit regime $n\sim N^2$.

Assign to $\la\in\Y_n(N)$ an $N$--particle configuration on $\Z_+$:
$$
\XX(\la)=\{x_1,\dots,x_N\}, \qquad x_i=\la_i+N-i  \tag2.2
$$
(note a difference from \tht{1.1}; due to the restriction $\ell(\la)\le N$,
$\la$ is uniquely determined by $\XX(\la)$). Then the probability space
$(\Y_n(N),M^\SW_{n,N})$ gives rise to an ensemble of random $N$--particle
configurations on $\Z_+$. More generally, we will deal with ensembles of
infinite random particle configurations as well. Such ensembles are examples of
what is called a {\it random point process\/} (or random point field). For a
discrete state space $\X$ (in our concrete case $\X=\Z_+$), a random point
process in $\X$ is determined by specifying a probability measure $\Cal P$ on
the space $\Conf(\X)=2^{\X}$ of all subsets in $\X$. Note that $\Conf(\X)$ is a
compact topological space in the natural topology.\footnote{About random point
processes in general, see, e.g., \cite{DVJ}, \cite{So1}.}

\proclaim{Conjecture 2.1 (Biane)} Consider the random Young diagram $\la$
distributed according to the measure $M^\SW_{n,N}$ on $\Y_n(N)$. Assume that
$n\to\infty$ and
$$
N=n^{1/2}-\tfrac12s\cdot n^{1/4}+o(n^{1/4})
$$
where $s$ is an arbitrary fixed real number. Equivalently,
$$
n=N^2+(s+o(1))N^{3/2}.
$$

Then the random configuration $\XX(\,\cdot\,)$ converges to a nontrivial random
point process on $\Z_+$ depending on $s$.
\endproclaim

Here convergence means weak convergence of probability measures on the compact
space $\Conf(\Z_+)$. The limit process is nontrivial in the sense that the
limiting measure on $\Conf(\Z_+)$ does not reduce to the delta measure on the
empty configuration or on the configuration coinciding with the whole set
$\Z_+$.

In section 3 we introduce the random point processes that appear as limit
processes for Biane's model, and in section 4 we verify Biane's conjecture in a
modified formulation.

\head 3. The discrete Hermite kernel \endhead

 We start with some necessary
generalities. Let $\X$ be a discrete space and $\Cal P$ be a random point
process in $\X$ (that is, a probability measure on $\Conf(\X))$. The {\it
correlation functions\/} $\rho_n$ of $\Cal P$ ($n=1,2,\dots$) are probabilities
for random configurations $X\in\Conf(\X)$ to contain a given finite set
$\{x_1,\dots,x_n\}$:
$$
\rho_n(x_1,\dots,x_n)=\Prob\{\{x_1,\dots,x_n\}\subseteq X\}.
$$

The initial measure $\Cal P$ is uniquely determined by the correlation
functions $\rho_1,\rho_2,\dots$. Indeed, for any finite subset $A\subset\X$,
there is a natural projection $\Conf(\X)\to\Conf(A)$ given by taking
intersection: $X\mapsto X_A=X\cap A$. Let $\Cal P_A$ be the pushforward of
$\Cal P$ under this projection; this is a probability measure on the finite set
$\Conf(A)$. Using the inclusion--exclusion principle it is readily seen that
$\Cal P_A$ is determined by the values of the correlation functions on $A$. For
instance, for $A=\{a,b\}$ we have
$$
\gathered \Prob\{X_A=\{a,b\}\}=\rho_2(a,b)\\
\Prob\{X_A=\{a\}\}=\rho_1(a)-\rho_2(a,b),\\
\Prob\{X_A=\{b\}\}=\rho_1(b)-\rho_2(a,b),\\
\Prob\{X_A=\varnothing\}=1-\rho_1(a)-\rho_1(b)+\rho_2(a,b).
\endgathered
$$
On the other hand, the initial measure $\Cal P$ is clearly determined by
collection of the measures $\Cal P_A$.

We say that a sequence $\Cal P_1, \Cal P_2, \dots$ of random point processes in
$\X$ converges to a random point process $\Cal P$ in the same space, $\Cal
P_k\to\Cal P$, if the corresponding probability measures weakly converge. This
happens if and only if the correlation functions of the processes $\Cal P_k$
pointwise converge to the respective correlation functions of the process $\Cal
P$. Indeed, by the definition of the topology in $\Conf(\X)$, we have $\Cal
P_k\to\Cal P$ if and only $(\Cal P_k)_A\to\Cal P_A$ for any finite
$A\subset\X$, and the latter is clearly equivalent to convergence of
correlation functions.

A random point process in $\X$ is said to be {\it determinantal\/} if there
exists a function $K(x,y)$ on $\X\times\X$ such that the correlation functions
are given by the determinantal formula
$$
\rho_n(x_1,\dots,x_n)=\det[K(x_i,x_j)], \qquad n=1,2,\dots,
$$
where the determinant in right-hand side has order $n\times n$. Then $K$ is
called the {\it correlation kernel\/} of $\Cal P$. Thus, a determinantal
process is uniquely determined by its correlation kernel. If a kernel $K$ is
Hermitean--symmetric, $K(x,y)=\overline{K(y,x)}$, then it serves as a
correlation kernel of a random point process if and only if $\Vert K\Vert\le1$,
that is, $K$ corresponds to a contractive operator in the Hilbert space
$\ell^2(\X)$. Indeed, this is a very special case of \cite{So1, Thm. 3}.

In particular, any projection kernel (that is, the kernel corresponding to a
selfadjoint projection operator in $\ell^2(\X)$) determines a random point
process in $\X$. We will introduce now a concrete family of projection kernels
which we will need in the sequel.

Consider the semi--infinite Jacobi matrix
$$
D^\he=\bmatrix 0 & \sqrt1 & 0 & 0 & \hdots\\
\sqrt1 & 0 & \sqrt2 & 0 & \hdots\\
0 & \sqrt2 & 0 & \sqrt 3 & \hdots\\
0 & 0 & \sqrt 3 & 0 & \dots\\
\vdots & \vdots & \vdots & \vdots & \ddots
\endbmatrix \tag3.1
$$
(the origin of this matrix and of its notation will become clear soon). We
agree that the rows and columns are indexed by the nonnegative integers
$x\in\Z_+$. Then the matrix $D^\he$ determines a symmetric operator in the
Hilbert space $\ell^2(\Z_+)$: by definition, the domain of the operator is the
algebraic span of the basis elements $\{\de_x\}$, $x\in\Z_+$. We will denote
this operator by the same symbol $D^\he$.

\proclaim{Lemma 3.1} The operator $D^\he$ is essentially self--adjoint. Its
closure $\overline{D^\he}$ has simple purely continuous Lebesgue spectrum. For
any Borel set $B\subseteq\R$, the corresponding spectral projection operator
$P_{B}$ is given by the kernel
$$
P_B(x,y)=(P_B\de_y,\de_x)=(2\pi x!y!2^{x+y})^{-1/2}\, \int_B
e^{-t^2/2}H_x(t/\sqrt2)H_y(t/\sqrt2)d\si, \tag3.2
$$
where $x,y$ range over $\Z_+$ and $H_x(t)$ denotes the Hermite polynomial of
degree $x$.
\endproclaim

\demo{Proof} We will need a few basic facts concerning the classical moment
problem:

Let $\rho$ be a measure on $\R$ with infinite support and with finite moments
of all orders. Then the space $\C[t]$ of polynomials in one variable can be
viewed as a subspace of the Hilbert space $L^2(\R,\rho)$. Let
$\overline{\C[t]}$ denote the closure of this subspace. Finally, let $m_n=\int
t^m \rho(dt)$ be the moments of $\rho$, $n=0,1,2,\dots$.

(A) The operator of multiplication by $t$ with domain of definition $\C[t]$ is
essentially self--adjoint in $\overline{\C[t]}$ if and only if the moment
problem associated with the sequence $\{m_n\}$ is determinate (that is, $\rho$
is a unique measure on $\R$ with moments $m_n$). See, e.g., \cite{Si, p.~86,
Thm. 2}.

(B) If the moment problem associated with $\{m_n\}$ is determinate then
$\overline{\C[t]}$ coincides with the whole space $L^2(\R,\rho)$. See, e.g.,
\cite{Si, p.~131, Prop.~4.15} or \cite{Ak, Cor.~2.3.3}.

(C) The moment problem associated with $\{m_n\}$ is determinate if the moments
grow not too fast. For instance, a simple sufficient condition says that the
moment problem is determinate if the exponential generating series
$$
g(z)=\sum_{n=0}^\infty m_n\frac{z^n}{n!}
$$
is analytic in a neighborhood of $z=0$ (which holds if the function $t\mapsto
e^{zt}$ is $\rho$--integrable for sufficiently small $z$). See, e.g., \cite{Si,
p.~88, Prop.~1.5}.

Finally, let $p_0,p_1,\dots$ stand for the orthogonal polynomials with respect
to $\rho$, normalized so that $\int p^2_n(t)\rho(dt)=1$. Then $\{p_n\}$ is an
orthonormal basis in $\overline{\C[t]}$. The polynomials $p_n$ satisfy a
three--term recurrence relation, which means that in the basis $\{p_n\}$, the
matrix of the operator of multiplication by $t$ is a (symmetric) tridiagonal
matrix.

Now we return to the proof of the lemma. Take as $\rho$ the normal distribution
$$
\rho^\he(dt)=\frac1{\sqrt{2\pi}}\,e^{-t^2/2}dt.
$$
The corresponding polynomials $p_n$ are
$$
\wt H_n(t)=(n!\,2^n)^{-1/2}H_n(t/\sqrt 2),
$$
where the $H_n$'s are the Hermite polynomials in the standard normalization,
see \cite{KS, \S1.13}. The three--term recurrence relation for the $H_n$'s has
the form
$$
H_{n+1}(t)-2tH_n(t)+2nH_{n-1}(t)=0, \tag3.2
$$
see \cite{KS, \S1.13}. Rewriting this in terms of the $\wt H_n$'s we get
$$
t\wt H_n(t)=\sqrt{n+1}\wt H_{n+1}(t)+\sqrt{n}\wt H_{n-1}(t).
$$
Thus, the matrix of multiplication by $t$ in the basis $\{\wt H_n\}$ is just
the Jacobi matrix $D^\he$ as defined in \tht{3.1}.

It is readily checked that $\rho^\he$ satisfies condition (C) above. By (B),
the space of polynomials is dense in $L^2(\R,\rho^\he)$. Consequently, $\{\wt
H_n\}$ is an orthonormal basis in $L^2(\R,\rho^\he)$, and the correspondence
$\de_x\leftrightarrow \wt H_x$ makes it possible to identify the Hilbert spaces
$\ell^2(\Z_+)$ and $L^2(\R,\rho^\he)$. By (A), our operator $D^\he$ is
essentially self--adjoint. \footnote{Of course, all these facts are well
known.} Now it is clear that we have obtained the explicit spectral
decomposition for the corresponding self--adjoint operator $\overline{D^\he}$,
the closure of $D^\he$. Namely, for any bounded real Borel function $\chi$ on
$\R$, the operator $\chi(\overline{D^\he})$ is realized as the operator of
multiplication by $\chi$ in  $L^2(\R,\rho^\he)$. In particular, taking as
$\chi$ the characteristic function $\chi_B$ of a Borel set $B\subset\R$ we see
that the corresponding spectral projection $P_B$ is the operator of
multiplication by $\chi_B$. The matrix of $P_B$ in the basis $\{\wt H_n\}$ is
given by formula \tht{3.2}, which concludes the proof. \qed
\enddemo

\example{Definition 3.2} The {\it discrete Hermite kernel\/} with parameter
$s\in\R$ is the above spectral kernel corresponding to the set $B=[s,\infty)$.
That is
$$
K^\he_s(x,y)=(2\pi x!y!2^{x+y})^{-1/2}\,
\int_s^{+\infty}e^{-t^2/2}H_x(t/\sqrt2)H_y(t/\sqrt2)dt, \tag3.3
$$
where $x,y\in\Z_+$.
\endexample

The next proposition provides alternative expressions for this kernel:

\proclaim{Proposition 3.3} For $x\ne y$, the discrete Hermite kernel can also
be written as
$$
(\pi x!y!2^{x+y})^{-1/2}\,e^{-s^2/2} \cdot \frac{xH_{x-1}(s/\sqrt2)\cdot
H_y(s/\sqrt2) -H_x(s/\sqrt2)\cdot yH_{y-1}(s/\sqrt2)}{x-y} \tag3.4
$$
or equivalently as
$$
(4\pi x!y!2^{x+y})^{-1/2}\,e^{-s^2/2} \cdot \frac{H_{x+1}(s/\sqrt2)\cdot
H_y(s/\sqrt2) -H_x(s/\sqrt2)\cdot H_{y+1}(s/\sqrt2)}{x-y} \tag3.5
$$
\endproclaim

\demo{Proof} The equivalence of \tht{3.4} and \tht{3.5} follows from the
three--term relation \tht{3.2}.

Next, we will employ two relations for Hermite polynomials (see \cite{KS,
\tht{1.13.6} and \tht{1.13.8}}):
$$
\gather (H_{n+1}(t))'=2(n+1)H_n(t), \tag3.6\\
(e^{-t^2}H_n(t))'=-e^{-t^2}H_{n+1}(t). \tag3.7
\endgather
$$

Abbreviating
$$
C(x,y)=(2\pi x!y!2^{x+y})^{-1/2}
$$
we have
$$
\aligned K^\he_s(x,y)&=C(x,y)
\int_s^{+\infty}e^{-t^2/2}H_x(t/\sqrt2)H_y(t/\sqrt2)dt\\
&=C(x,y)\sqrt2 \int_{s/\sqrt2}^{+\infty}e^{-t^2}H_x(t)H_y(t)dt
\endaligned
$$
Multiplying by $x-y=(x+1)-(y+1)$ and using \tht{3.6} we get
$$
(x-y)K^\he_s(x,y)=\frac{C(x,y)}{\sqrt2}
\int_{s/\sqrt2}^{+\infty}e^{-t^2}\left(H'_{x+1}(t)H_y(t)-H_x(t)H'_{y+1}(t)\right)dt.
$$
Integrating by parts and using \tht{3.7} we finally get
$$
\multline
(x-y)K^\he_s(x,y)\\=\frac{C(x,y)}{\sqrt2}e^{-s^2/2}\left(H_{x+1}(s/\sqrt2)H_y(s/\sqrt2)
-H_x(s/\sqrt2)H_{y+1}(s/\sqrt2)\right),
\endmultline
$$
which equals the expression \tht{3.5} multiplied by $x-y$. \qed
\enddemo

\head 4. Proof of modified Biane's conjecture \endhead

We will apply a well--known trick called {\it poissonization\/}. Its general
idea is to make a large parameter $n$ random and obeying the Poisson
distribution on $\Z_+$ with large parameter $\nu$. Due to the asymptotic
concentration of the Poisson distribution one believes that the large $n$ limit
regime and the large $\nu$ limit regime are equivalent (of course, this claim
has to be justified in each concrete situation). On the other hand, the latter
regime often turns out to be easier to study.

For instance, as shown in \cite{Jo1} and \cite{BOO}, application of the
poissonization procedure to the Plancherel measures $M^\Plancherel_n$ leads to
determinantal point processes. The same happens for the measures $M^\SW_{n,N}$
(Lemma 4.2 below).

By definition, the poissonized version of $M^\SW_{n,N}$, denoted as
$M^\PSW_{\nu,N}$, lives on the set
$$
\Y(N)=\bigcup_{n=0}^\infty \Y_n(N)=\{\la\in\Y\mid \ell(\la)\le N\}
$$
of all Young diagrams with at most $N$ rows. This new measure depends on a
positive parameter $\nu$, and is given by
$$
M^\PSW_{\nu,N}(\la)=e^{-\nu}\,\frac{\nu^{|\la|}}{|\la|!}\,
M^\SW_{|\la|,N}(\la), \qquad \la\in\Y(N). \tag4.1
$$
Clearly, $M^\PSW_{\nu,N}$ is a probability measure. In the present paper we do
not justify the poissonization procedure and simply replace the measures
$M^\SW_{n,N}$ by their poissonized versions in Biane's conjecture. Theorem 4.1
stated below proves the conjecture and identifies the limit process.

Let $X_{\nu,N}$ be the random $N$--particle configuration on $\Z_+$ obtained
via the correspondence \tht{2.2} from the random Young diagram $\la$
distributed according to the measure $M^\PSW_{\nu,N}$ on $\Y(N)$. That is, if
$\{x_1,\dots,x_N\}=\wt X(\la)$ then
$$
\Prob(\{x_1,\dots,x_N\})=M^\PSW_{\nu,N}(\la). \tag4.2
$$

\proclaim{Theorem 4.1} Fix an arbitrary $s\in\R$. Let $N=1,2,\dots$ and assume
that the parameter $\nu=\nu(N)$ depends on $N$ in such a way that
$$
\nu=N^2+(s+o(1))N^{3/2}, \qquad N\to\infty. \tag4.3
$$

As $N\to\infty$, $X_{\nu(N),N}$ converges to the determinantal point process on
$\Z_+$ with the correlation kernel $K^\he_s(x,y)$ as defined in\/ {\rm \S3}.
\endproclaim

\demo{Proof} It suffices to verify that the correlation functions of
$X_{\nu(N),N}$ converge to the respective correlation functions given by
correlation kernel $K^\he_s(x,y)$ (see the beginning of \S3). To do this we
will prove that $X_{\nu(N),N}$ is a determinantal process (Lemma 4.2 below) and
its correlation kernel pointwise converges to $K^\he_s(x,y)$ (Lemma 4.4 below),
which implies the claim of the theorem.

Consider the weight function for the Charlier polynomials with parameter
$\th>0$:
$$
W_\th^\ch(x)=\frac{\th^x}{x!}\,, \qquad x\in\Z_+\,,
$$
see \cite{KS, \S1.12}. The $N$--particle {\it Charlier ensemble\/} is formed by
random $N$--particle configurations $\{x_1,\dots,x_N\}\subset\Z_+$ such that
$$
\Prob(\{x_1,\dots,x_N\})=\const(N,\th)\cdot\prod_{i=1}^N W_\th^\ch(x_i)\cdot
\prod_{1\le i<j\le N}(x_i-x_j)^2, \tag4.4
$$
where $\const(N,\th)$ is the normalization constant (it can be evaluated
explicitly but we do not need the precise expression).

Let $C_m(x;\th)$ denote the Charlier polynomial of degree $m$, and let $\Vert
C_m(\,\cdot\,;\th)\Vert$ be its norm in the weighted $\ell^2$ space with the
weight function $W_\th^\ch$:
$$
\Vert C_m(\,\cdot\,;\th)\Vert^2=\sum_{x=0}^\infty C_m^2(x;\th)W_\th^\ch(x).
$$
The normalized functions
$$
\wt C_m(x;\th)=(W^\ch(x))^{1/2}\Vert C_m(\,\cdot\,;\th)\Vert^{-1}\,C_m(x;\th),
\qquad m=0,1,2,\dots
$$
form an orthonormal system in the ordinary space $\ell^2(\Z_+)$.

As is well known, the Charlier ensemble is a determinantal point process with
the correlation kernel
$$
K^\ch_{N,\th}(x,y)=\sum_{m=0}^{N-1}\wt C_m(x;\th)\wt C_m(y;\th).
$$
See, e.g., \cite{Jo1}, \cite{K\"o}. This is a projection kernel: the
corresponding operator is the projection in $\ell^2(\Z_+)$ on the
$N$--dimensional subspace spanned by the functions $\wt C_0, \dots,\wt
C_{N-1}$.

\proclaim{Lemma 4.2} For any $\nu>0$ and $N=1,2,\dots$, the random process
$X_{\nu,N}$ coincides with the $N$--particle Charlier ensemble with parameter
$\th=\nu/N$.
\endproclaim

\demo{Proof} Let us compare the right--hand sides of \tht{4.2} and \tht{4.4}.
The right--hand side of \tht{4.2} is defined by \tht{2.1} and \tht{4.1}, this
gives
$$
e^{-\nu}\,\frac{\nu^{|\la|}\dim\la\,\Dim_N \la}{|\la|!N^{|\la|}}.
$$
We have
$$
\frac{\dim\la}{|\la|!}=\frac{\prod\limits_{1\le i<j\le
N}(x_i-x_j)}{\prod\limits_{1\le i\le N}x_i!} \tag4.5
$$
(Frobenius' formula, see, e.g., \cite{Ma}) and
$$
\Dim_N\la=\prod_{1\le i<j\le N}\frac{x_i-x_j}{j-i}
$$
(Weyl's character formula). Finally,
$$
|\la|=x_1+\dots+x_N-\tfrac{N(N-1)}2.
$$
Using these expressions we obtain the right--hand side of \tht{4.4} with an
appropriate constant. \qed
\enddemo

Given $\th>0$, consider the semi--infinite Jacobi matrix
$$
D_\th^\ch=\bmatrix 0 & \sqrt1 & 0 & 0 & \hdots\\
\sqrt1 & -1/\sqrt\th & \sqrt2 & 0 & \hdots\\
0 & \sqrt2 & -2/\sqrt\th & \sqrt 3 & \hdots\\
0 & 0 & \sqrt 3 & -3/\sqrt\th & \dots\\
\vdots & \vdots & \vdots & \vdots & \ddots
\endbmatrix
$$
We also regard it as a symmetric operator in $\ell^2(\Z_+)$ whose domain of
definition is the set of finite linear combination of the basis elements.

\proclaim{Lemma 4.3} The operator $D_\th^\ch$ is essentially self--adjoint. Its
closure $\overline{D_\th^\ch}$ has purely point spectrum
$\{\frac{\th-m}{\sqrt\th}\mid m=0,1,\dots\}$. The kernel $K^\ch_{N,\th}$
coincides with the kernel of the spectral projection operator corresponding to
the part of spectrum
$$
\left\{\frac{\th-m}{\sqrt\th}\,, \quad m=0,1,\dots, N-1\right\}. \tag4.6
$$
\endproclaim

\demo{Proof} The Charlier polynomials with parameter $\th$ satisfy the
difference equation
$$
\th C_m(x+1;\th)-x C_m(x;\th)+ x C_m(x-1;\th)=(\th-m) C_m(x;\th),
$$
see \cite{KS, \tht{1.12.5}}. {}From the expression for the weight function it
follows that for the normalized functions $\wt C_m(x;\th)$, this equation is
transformed into
$$
\sqrt{x+1}\wt C_m(x+1;\th)-\frac{x}{\sqrt\th}\wt C_m(x;\th)+\sqrt x\wt
C_m(x-1;\th)=\frac{\th-m}{\sqrt\th}\wt C_m(x;\th).
$$
We see that $D_\th^\ch$ is precisely the difference operator on $\Z_+$ standing
in the left--hand side.

We claim that $D_\th^\ch$ is essentially self--adjoint and $\{\wt C_m\}$ is the
complete set of the eigenvectors of the self--adjoint operator
$\overline{D_\th^\ch}$. Indeed, the three--term recurrence relation
$$
-xC_m(x;\th)=\th C_{m+1}(x;\th)-(m+\th)C_m(x;\th)+mC_{m-1}(x;\th)
$$
(see \cite{KS, \tht{1.12.3}}) and the explicit expression for the norm
$$
\Vert C_m\Vert^2=\th^{-m}e^\th m!
$$
(see \cite{KS, \tht{1.12.2}}) imply that the same Jacobi matrix corresponds to
the operator of multiplication by $(\th-x)/\sqrt\th$ in the basis $\{C_m/\Vert
C_m\Vert \}$ of the space of polynomials. It is readily checked that the
Charlier weight, viewed as a measure on $\Z_+$, satisfies the sufficient
condition (C), see the proof of Lemma 3.1. Then the same argument as in that
lemma shows that the space of polynomials in dense in the weighted space
$\ell^2(\Z_+,W_\th^\ch)$ and the above multiplication operator is essentially
self--adjoint.

This is equivalent to saying that the functions $\wt C_m$ form an orthonormal
basis in $\ell^2(\Z_+)$ and $D_\th^\ch$ is essentially self--adjoint. Then it
follows from the difference equation that $\overline{D_\th^\ch}$ has $\wt C_m$
as an eigenvector with eigenvalue $\frac{\th-m}{\sqrt\th}$. The last claim of
the lemma is now obvious. \qed
\enddemo

\proclaim{Lemma 4.4} Let $s\in\R$ be fixed and assume
$\th=\th(N)=N+(s+o(1))N^{1/2}$. Then
$$
\lim_{N\to\infty} K^\ch_{N,\th(N)}(x,y)=K^\he_s(x,y), \qquad x,y\in\Z_+\,.
$$
\endproclaim

Here the assumption on $\th$ comes from \tht{4.3} and the relation $\th=\nu/N$
(Lemma 4.2).

\demo{Proof} Consider the self--adjoint operators $\overline{D_\th^\ch}$ (where
$\th>0$) and $\overline{D^\he}$ in $\ell^2(\Z_+)$. Let $\ell^2_0(\Z_+)$ denote
the algebraic linear span of the basis elements $\de_x$, $x\in\Z_+$. By Lemma
3.1 and Lemma 4.3, all these operators have $\ell^2_0(\Z_+)$ as a common
essential domain. Moreover, it is evident that if $\th\to\infty$ then
$D_\th^\ch\to D^\he$ on $\ell^2_0(\Z_+)$. It follows that
$\overline{D_\th^\ch}\to\overline{D^\he}$ in the strong resolvent sense (see
\cite{RS, Thm. VIII.25}).

Let us regard $K^\ch_{N,\th}$ and $K^\he_s$ as operators in $\ell^2(\Z_+)$. By
Lemma 3.1, the latter operator is the spectral projection of $\overline{D^\he}$
corresponding to the set $[s,+\infty)$. By Lemma 4.3,  the former operator is
the spectral operator of $\overline{D_\th^\ch}$ corresponding to the set
\tht{4.6}. Next, it follows from the description of the spectrum of
$\overline{D_\th^\ch}$ in Lemma 4.3 that instead of the finite set \tht{4.6} we
can equally well take the continuous interval
$$
\left[\frac{\th-N+1}{\sqrt\th}\,, \quad +\infty\right).
$$

If $\th=\th(N)=N+(s+o(1))N^{1/2}$ then the left end of this interval can be
written as $s+\varepsilon_N$ where $\varepsilon_N\to0$ as $N\to\infty$. Since
$\overline{D^\he}$ has purely continuous spectrum, the strong resolvent
convergence implies that the spectral projection of $\overline{D_{\th(N)}^\ch}$
corresponding to $[s+\varepsilon_N,+\infty)$ strongly converges to the spectral
projection of $\overline{D^\he}$ corresponding to $[s,+\infty)$: this is proved
exactly as claim (b) in \cite{RS, Thm. VIII.24}. \qed
\enddemo

Note that Lemma 4.4 could be obtained from the known asymptotics for the
Laguerre polynomials \cite{Te} and the well--known connection between the
Laguerre and Charlier polynomials.

Lemma 4.4 completes the proof of Theorem 4.1. \qed
\enddemo

\head 5. Another model \endhead

Let $N$ and $M$ be two natural numbers, and $(M^N)$ be the rectangular Young
diagram with $N$ rows and $M$ columns. Given a Young diagram
$\la\subseteq(M^N)$, we denote by $(M^N)/\la$ the skew diagram which is the
difference of $(M^N)$ and $\la$. Reading this skew diagram from the bottom to
the top we get an ordinary Young diagram which will be denoted by
$\widehat\la$:
$$
(\widehat\la_1,\dots,\widehat\la_N)=(M-\la_N,\dots,M-\la_1).
$$

Let $\pi_{(M^N)}$ denote the irreducible representation of the symmetric group
$S_{NM}$ of degree $NM$, indexed by $(M^N)$. For any  $n=0,1,\dots,NM$, the
restriction of $\pi_{(M^N)}$ to the Young subgroup $S_n\times S_{NM-n}$ has
simple spectrum consisting of the irreducible representations of the form
$\pi_{\la}\otimes\pi_{\widehat\la}$ (the outer tensor product of the
irreducible representations of $S_n$ and $S_{NM-n}$, indexed by $\la$ and
$\widehat\la$, respectively). Indeed, this follows from the fact that the skew
Schur function $s_{(M^N)/\la}$ equals the ordinary Schur function
$s_{\widehat\la}$\,, as is readily verified using the Jacobi--Trudi formula
(see \cite{Ma, Ch. I, \tht{5.4}}).

Let $\Y_n(N,M)$ stand for the set of Young diagrams with $n$ squares, contained
in the rectangle $(M^N)$, $n=0,1,\dots,NM$. The above claim shows that the
following expression defines a probability measure on $\Y_n(N,M)$, which will
be denoted as $M_{n,N,M}$:
$$
M_{n,N,M}(\la)=\frac{\dim\la\cdot\dim\widehat\la}{\dim(M^N)}\,, \qquad
\la\in\Y_n(N,M).
$$

It turns out that if the triple of parameters $n,N,M$ goes to infinity in an
appropriate way then the boundary of the random Young diagram distributed
according to the measure $M_{n,N,M}$, after a suitable scaling, tends to a
nonrandom curve: This is a particular case of the results in Biane \cite{Bi1,
Thm. 3.1.2} and Pittel and Romik \cite{PR, \S1.1 and \S1.5}. Biane's approach
uses free probability. The derivation of Pittel and Romik of the explicit form
of the limit curve is based on the variational principle. Here we sketch a
simple alternative argument. It does not rigorously prove the existence of the
limit curve but allows one to guess what it is.

In what follows we will assume $M=N$ and abbreviate $M_{n,N}=M_{n,N,N}$. The
case of a rectangle $(M^N)$ can be handled in a similar way. We stick to the
square case $M=N$ to simplify the notation only.

Assume that $N$ and $n$ go to infinity in such a way that $n\sim pN^2$, where
$p\in(0,1)$ is a fixed parameter. Instead of $M_{n,N}$ we will be dealing with
a modified measure, which is obtained by a mixing procedure similar to
poissonization: all values $n=0, 1,\dots,N^2$ are mixed by making use of the
binomial distribution on $\{0,1,\dots,N^2\}$ with parameter $p$. Like the
Poisson distribution, the binomial distribution possesses the concentration
property: as $N$ gets large, the main contribution comes from those $n$'s which
are close to $pN^2$. Thus, one may believe that mixing does not affect the
asymptotics.

The resulting measure lives on the set $\Y(N,N)$ of all Young diagrams
contained in $(N^N)$ (no constraints on $|\la|$ are imposed):
$$
M^{\Mix}_{p,N}(\la) =\binom{N^2}{|\la|}p^{|\la|}(1-p)^{N^2-|\la|}\, M_{|\la|,
N}(\la), \qquad \la\in\Y(N,N). \tag5.1
$$
The next claim, which  is similar to Lemma 4.2, shows that the measure
$M^{\Mix}_{p,N}$ leads to the $N$--particle {\it Krawtchouk ensemble\/}.

Denote by $W^\kr_{p,L}$ the weight function of the Krawtchouk orthogonal
polynomials on the finite set of integers $\{0,1,\dots,L\}$ and depending on
the parameter $p\in(0,1)$:
$$
W^\kr_{p,L}(x)=\binom{L}{x}p^x(1-p)^{L-x}, \qquad x=0,1,\dots,L.
$$

\proclaim{Lemma 5.1} Under the correspondence $\la\leftrightarrow
\{x_1,\dots,x_N\}$ defined by \tht{2.2}, random Young diagrams $\la\in\Y(N,N)$
distributed according to the measure $M^{\Mix}_{p,N}$ turn into random
$N$--particle configurations in $\{0,1,\dots,L\}$ with $L=2N-1$ and such that
$$
\Prob(\{x_1,\dots,x_N\})=\const(p,N)\cdot\prod_{i=1}^N W^\kr_{p,L}(x_i)
\prod_{1\le i<j\le N}(x_i-x_j)^2. \tag5.2
$$
\endproclaim

\demo{Proof} Recall that
$$
M_{|\la|, N}(\la)=\frac{\dim\la\cdot\dim(\widehat\la)}{\dim(N^N)}\,.
$$
Applying Frobenius' dimension formula \tht{4.5} and using \tht{5.1} we obtain
the desired expression. \qed
\enddemo

The next result describes the limit behavior of the Krawtchouk ensemble ``in
the bulk''.

Fix a real number $c$ such that $|c|<2\sqrt{p(1-p)}$.

Let $a_N$ be an arbitrary sequence of integers such that $a_N\sim cN$. Given a
random configuration $\{x_1,\dots,x_N\}$ of the $N$--particle Krawtchouk
ensemble \tht{5.2}, we shift it by $N+a_N$ to the left and obtain in this way a
new random $N$--particle configuration $\{x'_1,\dots,x'_N\}\subset\Z$:
$$
x'_i=x_i-N-a_N, \qquad 1\le i\le N. \tag5.3
$$

\proclaim{Proposition 5.2} Under the above assumptions,  the random
configuration \tht{5.3} converges as $N\to\infty$ to the translation invariant
determinantal random point process on $\Z$ with the correlation kernel
$$
K^\dsine_\varphi(x,y)=\frac{\sin(\varphi(x-y))}{\pi(x-y)}\,, \qquad
x,y\in\Z,\tag5.4
$$
where
$$
\varphi=\arccos\left(\frac{c(1-2p)}{2\sqrt{(1-c^2)p(1-p)}}\right). \tag5.5
$$
\endproclaim

The assumption $|c|<2\sqrt{p(1-p)}$ just means that $
\left|\frac{c(1-2p)}{2\sqrt{(1-c^2)p(1-p)}}\right|<1$, so that $\varphi$ is
well defined.

The kernel \tht{5.4}, called the {\it discrete sine kernel\/}, first appeared
in connection with the Plancherel model, see \cite{BOO}. This kernel is a
lattice counterpart of the celebrated sine kernel
$$
K^{\text{sine}}(u,v)=\frac{\sin(\pi(u-v))}{\pi(u-v)}\,, \qquad u,v\in\R.
$$
The result of the proposition is a manifestation of a general phenomenon
studied in \cite{BKMM}: The discrete sine kernel is the universal correlation
kernel arising in the bulk limit of discrete orthogonal polynomial ensembles.

\demo{Sketch of proof} The argument is very similar to the proof of Theorem
4.1. We give the formal computation below; the justification is omitted.

Let $K_m(x;p,L)$ denote the Krawtchouk polynomial of degree $m$ (necessary
information about these polynomials can be found in \cite{KS, \S1.10}). The
normalized functions
$$
\wt K_m(x;p,L)=(W^\kr_{p,L}(x))^{1/2}\,\Vert K_m(\,\cdot\,;p,L)\Vert^{-1}\,
K_m(x;p,L)
$$
form an orthonormal basis in the $\ell^2$ space on the finite set
$\{0,1,\dots,L\}$. The $N$--particle Krawtchouk ensemble \tht{5.2} is a
determinantal point process with the correlation kernel
$$
K^\kr_{p,L}(x,y)=\sum_{m=0}^{N-1}\wt K_m(x;p,L)\wt K_m(y;p,L). \tag5.6
$$

The Krawtchouk polynomials $K_m(x;p,L)$ satisfy the difference equation
$$
\multline
p(L-x)K_m(x+1;p,L)+x(2p-1)K_m(x;p,L)+x(1-p)K_m(x-1;p,L)\\
=(pL-m)K_m(x;p,L).
\endmultline
$$
In terms of the normalized functions the difference equation takes the form
$$
\multline \frac{\sqrt{(L-x)(x+1)}}L\,\wt K_m(x+1;p,L)
+\frac{x(2p-1)}{L\sqrt{p(1-p)}}\,\wt K_m(x;p,L)\\
+\frac{\sqrt{(L-x+1)x}}L\,\wt K_m(x-1;p,L) =\frac{pL-m}{L\sqrt{p(1-p)}} \wt
K_m(x;p,L).
\endmultline
$$

Let $D$ denote the difference operator defined by the left--hand side of this
equation. The correlation kernel \tht{5.6} corresponds to the projection on the
following part of the spectrum of the operator $D$:
$$
\left\{\frac{pL-m}{L\sqrt{p(1-p)}}\,,\qquad m=0,\dots,N-1\right\}. \tag5.7
$$

Recall that $L=2N-1$ and $x\approx N+cN+x'$. For large $N$, the three
coefficients of our difference operator are approximately equal to
$$
\frac12\sqrt{1-c^2}\,, \qquad \frac{(1+c)(2p-1)}{2\sqrt{p(1-p)}}\,, \qquad
\frac12\sqrt{1-c^2}\,,
$$
and the set \tht{5.7} approximates the interval
$$
\left[\frac{2p-1}{2\sqrt{p(1-p)}}\,, \qquad \frac{2p}{2\sqrt{p(1-p)}}\right].
\tag5.8
$$
Thus, in the limit $N\to\infty$ we get the difference operator
$$
\frac12\sqrt{1-c^2}f(x'+1)+\frac{(1+c)(2p-1)}{2\sqrt{p(1-p)}}f(x')+
\frac12\sqrt{1-c^2}f(x'-1), \qquad x'\in\Z,
$$
and the spectral projection corresponding to the interval \tht{5.8}.

Subtracting the scalar operator $f\mapsto \frac{(1+c)(2p-1)}{2\sqrt{p(1-p)}}f$
and dividing by $\frac12\sqrt{1-c^2}$ we finally arrive to the difference
operator
$$
f(x'+1)+f(x'-1), \qquad x'\in\Z, \tag5.9
$$
and the spectral interval
$$
\left[\frac{c(1-2p)}{\sqrt{(1-c^2)p(1-p)}}\,, \qquad
\frac{c(1-2p)+1}{\sqrt{(1-c^2)p(1-p)}}\right]. \tag5.10
$$

The corresponding spectral projection is given by the discrete sine kernel
\tht{5.4}. Indeed, to study the difference operator \tht{5.9} it is convenient
to make the Fourier transform from $\ell^2(\Z)$ to the $L^2$ space on the unit
circle $|z|=1$. Then \tht{5.9} becomes the operator of multiplication by the
function $z+\bar z=2\Re z$. Hence, we see that our operator has purely
continuous (double) spectrum ranging from $-2$ to 2. It is readily seen that
the right end of the interval \tht{5.10} is always $\ge2$, while the left end
is somewhere inside this interval (here we use the assumption
$|c|<2\sqrt{p(1-p)}$\,). Thus, in the $L^2$ space on the circle, our spectral
projection becomes the operator of multiplication by the characteristic
function of the arc going in the counterclockwise direction from
$e^{-i\varphi}$ to $e^{i\varphi}$, where $\varphi$ is given by \tht{5.5}. In
the $\ell^2(\Z)$--realization, this is the integral operator with the discrete
sine kernel \tht{5.4}.\qed
\enddemo

\proclaim{Corollary 5.3} Let $\la\in\Y(N,N)$ be the random Young diagram
distributed according to the probability measure $M^{\Mix}_{p,N}$. Assuming
that the boundary of $\la$ in the scaled coordinates $x=i/N$, $y=j/N$ has a
nonrandom limit described by a curve $x+y=F(y-x)$ we can explicitly find $F$
from the equation
$$
\frac{1-F'(c)}2=\frac{\varphi}\pi
=\frac1\pi\arccos\left(\frac{c(1-2p)}{2\sqrt{(1-c^2)p(1-p)}}\right),
$$
where $c$ ranges over the interval $(-2\sqrt{p(1-p)}, \quad2\sqrt{p(1-p)}\,)$.
\endproclaim

An additional condition is that the area bounded by the limit curve and the
coordinate axes in the $(x,y)$ plane has to be equal to $p$.

\demo{Idea of proof} We observe that the density function (that is, the first
correlation function) of the point process with discrete sine correlation
kernel is the constant $\varphi/\pi$. Then we use the same argument as in item
(b) of \S1 (see also \cite{BOO, Remark 1.7}). \qed
\enddemo

One can verify that this result agrees with the formulas in \cite{PR}. The
endpoints of the interval $(-2\sqrt{p(1-p)}, \quad2\sqrt{p(1-p)}\,)$
correspond to the endpoints of the limit curve that lie on the coordinate axes.

Note that the result of Proposition 5.2 can be obtained using asymptotics of
Krawtchouk polynomials obtained in \cite{IS}. The case $p=1/2$ is also handled
in \cite{Jo2, Lemma 2.8}. These papers contain much finer results on the
asymptotics but obtaining them requires substantially more work.

\bigskip

{\smc A.~Borodin}: Mathematics 253-37, Caltech, Pasadena, CA 91125, U.S.A.

E-mail address: {\tt borodin\@caltech.edu}

\bigskip

{\smc G.~Olshanski}: Dobrushin Mathematics Laboratory, Institute for
Information Transmission Problems, Bolshoy Karetny 19, 127994 Moscow GSP-4,
RUSSIA.

E-mail address: {\tt olsh\@online.ru}

\enddocument